\documentclass[9pt,a4paper]{amsart}
\usepackage{amsmath,amsthm,amssymb,color,latexsym, tikz, longtable, cite, float, subcaption, hyperref}
\usepackage[lmargin=30mm,rmargin=30mm,tmargin=35mm,bmargin=35mm]{geometry}
\usetikzlibrary {arrows.meta}

\newcommand{\spacing}[1]{
\renewcommand{\baselinestretch}{#1}
\setlength{\footnotesep}{\baselinestretch\footnotesep}}
\spacing{1.15}

\theoremstyle{plain}
\newtheorem{definition}{Definition}[section]
\newcommand{\floor}[1]{\lfloor#1\rfloor}
\newcommand{\ceil}[1]{\lceil#1\rceil}

\newcommand{\bc}[1]{\left\lceil#1\right\rceil}
\newcommand{\bfl}[1]{\left\lfloor#1\right\rfloor}
\newcommand{\vr}[1]{r(#1\mid W)}
\newcommand{\pr}[1]{\left(#1\right)}

\newcommand{\createNode}[2]{
    \node[circle, fill=red, opacity=0, scale=0.9] at #2 (#1_o) {};
    \node[circle, fill=lightgray, opacity=0.7,scale=0.5, draw=gray] at #2 (#1) {};
}

\newcounter{nodecount}

\newcommand{\drawpath}[4]{
    \foreach \i in {#1,...,#2}{
        \ifnum\i<#2
            \pgfmathtruncatemacro{\nexti}{\i+1}
            \draw[line width=#4 pt] (v#3\i) -- (v#3\nexti);
        \fi
    }
}

\newcommand{\createGrid}[2]{
    \foreach \i in {1,...,#1} {
        \foreach \j in {1,...,#1} {
            \pgfmathsetmacro\x{#2 + \i - 1}
            \pgfmathsetmacro\y{#2 + \j - 1}
            \createNode{c_\i_\j}{(\x, \y)}
        }
    }
}
\newcommand{\finiteGrid}[1]{
    \createGrid{#1}{1}
    \foreach \i in {1,...,#1} {
        \foreach \j in {1,...,#1} {
            \ifnum\i<#1
                \pgfmathtruncatemacro\nexti{\i+1}
                \draw[draw=gray, opacity=0.5] (c_\i_\j_o) -- (c_\nexti_\j_o);
            \fi
            \ifnum\j<#1
                \pgfmathtruncatemacro\nextj{\j+1}
                \draw[draw=gray, opacity=0.5] (c_\i_\j_o) -- (c_\i_\nextj_o);
            \fi
            \ifnum\i>1 \ifnum\j<#1
                \pgfmathtruncatemacro\previ{\i-1}
                \pgfmathtruncatemacro\nextj{\j+1}
                \draw[draw=gray, opacity=0.5] (c_\i_\j_o) -- (c_\previ_\nextj_o);
            \fi \fi

             \ifnum\i<#1\ifnum\j<#1
                \pgfmathtruncatemacro\nexti{\i+1}
                \pgfmathtruncatemacro\nextj{\j+1}
                \draw[draw=gray, opacity=0.5] (c_\i_\j_o) -- (c_\nexti_\nextj_o);
            \fi \fi
        }
    }
}

\newcommand{\infiniteGrid}{
    \createGrid{6}{1}
    \foreach \i in {1,...,6} {
        \foreach \j in {1,...,6} {
            \ifnum\i<6
            \ifnum \i= 4
            
            \else
                \pgfmathtruncatemacro\nexti{\i+1}
                \draw[draw=gray, opacity=0.5] (c_\i_\j_o) -- (c_\nexti_\j_o);
            \fi 
            \fi
            \ifnum\j<6 
            \ifnum \j= 2\else
                \pgfmathtruncatemacro\nextj{\j+1}
                \draw[draw=gray, opacity=0.5] (c_\i_\j_o) -- (c_\i_\nextj_o);
            \fi\fi
            \ifnum\i>1  \ifnum \j<6
            \ifnum \i = 5\else
            \ifnum \j = 2\else
                \pgfmathtruncatemacro\previ{\i-1}
                \pgfmathtruncatemacro\nextj{\j+1}
                \draw[draw=gray, opacity=0.5] (c_\i_\j_o) -- (c_\previ_\nextj_o);
            \fi \fi \fi \fi

             \ifnum\i<6 \ifnum\j<6 
             \ifnum \j = 2\else
             \ifnum \i= 4\else
                \pgfmathtruncatemacro\nexti{\i+1}
                \pgfmathtruncatemacro\nextj{\j+1}
                \draw[draw=gray, opacity=0.5] (c_\i_\j_o) -- (c_\nexti_\nextj_o);
            \fi \fi \fi\fi
        }
    }
    \foreach \i in {4} {
        \foreach \j in {1,...,6} {
            \pgfmathtruncatemacro\nexti{\i+1}
            \draw[dotted,draw=gray, opacity=0.5] (c_\i_\j_o) -- (c_\nexti_\j_o);
        }
    }
    \foreach \j in {2} {
        \foreach \i in {1,...,6} {
            \pgfmathtruncatemacro\nextj{\j+1}
            \draw[dotted,draw=gray, opacity=0.5] (c_\i_\j_o) -- (c_\i_\nextj_o);
        }
    }
}

\newcommand{\GraphGeneralizedThetaGraph}[1][]{%
\begin{tikzpicture}[scale=#1]
        \node[circle, fill=black, scale=0.5] at (0, 3) (c1) {};
        \node[above, yshift=3, scale = 0.7] at (c1) {$c_1$};
        \node[circle, fill=black, scale=0.5] at (0, -3) (c2) {};
        \node[below, yshift=-3, scale = 0.7] at (c2) {$c_2$};
        
        \draw[dotted,smooth] (-0.5,2) -- (0.5, 2);
        \draw[dotted,smooth] (-0.5,-2) -- (0.5, -2);
        \draw[dotted,smooth] (-0.5,0) -- (0.5, 0);
        
        \foreach \i in {0, 20, 40, 60, 70, 110, 120, 140, 160, 180} {
            \draw[domain=-3:-1.4, smooth, variable=\y] plot ({3*cos(\i)*sqrt(1-(\y/3)^2)}, \y);
            \draw[dotted, domain=-1.4:-0.4, smooth, variable=\y] plot ({3*cos(\i)*sqrt(1-(\y/3)^2)}, \y);
            \draw[dotted, domain=0.4:1.4, smooth, variable=\y] plot ({3*cos(\i)*sqrt(1-(\y/3)^2)}, \y);
            \draw[domain=1.4:3, smooth, variable=\y] plot ({3*cos(\i)*sqrt(1-(\y/3)^2)}, \y);
            
            \foreach \y in {-2.5, -2, -1.5, 1.5, 2, 2.5} {
                \node[circle, fill=black, scale=0.2] at ({3*cos(\i)*sqrt(1-(\y/3)^2)}, \y) {};
            }
        }
      \end{tikzpicture}
}

\newcommand{\PathExample}[1][]{%
    \begin{tikzpicture}
    \node[circle, fill=black, scale=0.5] (c1) at (0, 0) {};
    \node[above, yshift=3, scale=0.7] at (c1) {$c_1$};

    \node[circle, fill=black, scale=0.5] (v1) at (2, 0) {};
    \node[above, yshift=3, scale=0.7] at (v1) {$v_{i, 1}$};

    \node[circle, fill=black, scale=0.5] (v2) at (4, 0) {};
    \node[above, yshift=3, scale=0.7] at (v2) {$v_{i, 2}$};

    \node[circle, fill=black, scale=0.5] (v4) at (8, 0) {};
    \node[above, yshift=3, scale=0.7] at (v4) {$v_{i, s_i}$};

    \node[circle, fill=black, scale=0.5] (c2) at (10, 0) {};
    \node[above, yshift=3, scale=0.7] at (c2) {$c_2$};
    
    \draw (c1) -- (v1);
    \draw (v1) -- (v2);
    \draw[dotted] (4, 0) -- (5, 0);
    \draw[dotted] (7, 0) -- (8, 0);
    \draw (v4) -- (c2);
\end{tikzpicture}
}

\newtheorem{thm}{Theorem}[section]
\newtheorem{lem}[thm]{Lemma}
\newtheorem{cor}[thm]{Corollary}
\newtheorem{prop2}[thm]{Proposition}
\newtheorem{conj}[thm]{Conjecture}

\newtheorem{ob}[thm]{Observation}

\newcommand{\theorem}[1]{\begin{thm}#1\end{thm}}

\newcommand{\Proposition}[1]{\begin{prop2}#1\end{prop2}}
\newcommand{\proposition}[1]{\begin{prop2}#1\end{prop2}}

\newcommand{\lemma}[1]
{\begin{lem}#1\end{lem}}

\newcommand{\corollary}[1]{\begin{cor}#1\end{cor}}

\newcommand{\figref}[1]{Figure~\ref{fig:#1}}

\newcommand{\mmd}[2]{
#1\text{ MMD } #2}

\newcommand{\nmmd}[2]{
#1\text{ $\neg$MMD } #2}

\begin{document}

\title[Metric Dimension of Generalized Theta Graphs]{Metric Dimension of Generalized Theta Graphs}
\author[]{Nadia Benakli, Nicole Froitzheim, David Martinez}\

\thanks{This material is based upon work supported by the National Science Foundation under Grant No. DMS-2150251}

\keywords{theta graph, generalized theta graph, metric dimension}

\begin{abstract}
A vertex $w$ in a graph $G$ is said to resolve two vertices $u$ and $v$ if $d(w,u)\neq d(w, v)$.
A set $W$ of vertices is a \textit{resolving} set for $G$ if every pair of distinct vertices
is resolved by some vertex
in $W$. The \textit{metric dimension} of $G$ is the minimum
cardinality of such a set. In this paper, we investigate the metric dimension of generalized theta graphs, providing exact values and structural insights fir several subclasses.
\end{abstract}

\maketitle

\section{Introduction}
    The metric dimension, introduced by Slater \cite{Slater_Introduction} and independently by Harary and Melter \cite{Harary_Melter}, is a graph invariant used to distinguish vertices in a graph based on their distances to a fixed subset of vertices. It has practical applications in various disciplines, such as robot navigation, network localization, and chemical structure identification. This invariant, along with its many variants, has been extensively studied across several classes of graphs: for an overview, see the surveys \cite{Kuziak2021,Tillquist2023}.  In this paper, we consider a class of graphs called generalized theta graphs, which consists of two vertices connected by multiple internally disjoint paths. When there are exactly three such paths, the graph is called a theta graph. We denote a generalized theta graph by $\Theta(s_1, \ldots s_m)$, where $m$ is the number of internally disjoint paths and  $s_1$, \ldots $s_m$ denote the number of internal vertices on the $m$ paths, with the assumption that $s_1\leq s_2\leq \ldots \leq s_m$. Generalized theta graphs are a natural extension of cycles. The metric dimension of cycles is 2, and it has been determined for all unicyclic graphs \cite{Sedlar2022}. For base bicyclic graphs, those that do not contain vertices of degree 1 (i.e., leaves), the problem is categorized into three types. Type 1 bicyclic graphs consist of two cycles that share exactly one common vertex, and type 2 consist of two vertex-disjoint cycles connected by a path. The metric dimensions of types 1 and 2 have been determined \cite{Khan2023}. A bicyclic graph of type 3 is precisely a theta graph, denoted $\Theta(s_1,s_2,s_3)$, and its metric dimension has also been determined \cite{Knor2022remarks,Wang2024}. This is stated in the following theorem.
    \theorem{
Let $\Theta(s_1, s_2, s_3)$ be a theta graph. Then 
    \begin{align}
        \beta(\Theta(s_1, s_2, s_3)) = \begin{cases}
            3 &\, \text{if (1) $s_i=s$ or (2) $s_1=s_2$ and $s_3 = s_1 + 2$}\\
            2 &\, \text{otherwise}
        \end{cases}
    \end{align}
    
    }
T4
    
    This paper presents several results on the metric dimension of generalized theta graphs, $\Theta(s_1, \ldots s_m)$. In particular, we derive the lower and upper bounds given in the following theorem.
\theorem{
\label{thm:TotalBound}   
Given a generalized theta graph, $\Theta(s_1, s_2, \ldots, s_m)$, where $m \geq 2$, we have 
    \begin{align}
       m-3 \leq \beta(\Theta(s_1, s_2, \ldots, s_m)) \leq m
    \end{align}
}    

These bounds are not sharp, as will be shown in this paper. The paper is organized as follows. Section 2 is devoted to the basic definitions and preliminaries. In Section 3, we establish the bounds for the metric dimension of generalized theta graphs. In Section 4, we consider consecutive generalized theta graphs and derive results for particular subclasses, including the uniform case.  Additionally, since theta graphs (i.e., generalized theta graphs of multiplicity 3) and generalized theta graphs of multiplicity 4 are special, separate sections are devoted to these cases. In this paper, we provide a new proof of Theorem 1 and 2 for theta graphs \cite{Knor2022remarks}.

\section{Preliminaries}  
Throughout this paper, all graphs are assumed to be simple, undirected, connected, and finite. Given a graph $G$, we denote by $V(G)$ its set of vertices and by $E(G)$ its set of edges. When no confusion arises, we will simply use $V$ and $E$. The distance between two vertices $u$ and $v$ in $G$, denoted by $d(u,v)$, is the length of a shortest path connecting them, where the length of a path is the number of edges it contains. 

We also extend the concept of the distance function in a graph to consider a shortest path that traverses a sequence of vertices. Let $S= (v_1, v_2, \ldots, v_n)$ be a sequence of vertices in a graph $G$. We define $I[v_1, v_2, \ldots, v_n]$ to be the path obtained by concatenating, in order, shortest paths between each pair of consecutive vertices $v_i$, $v_{i+1}$. The length of this path is given by 
\begin{align}
    D(S)= \sum_{i=1}^{n-1} d(v_i, v_{i+1})
\end{align}
 Note that $D(S)$ need not equal the distance between the endpoints $v_1$ and $v_n$; in fact, $D(S) \geq  d(v_1, v_n)$. 

\begin{definition}
    Let $u,v$ and $w$ be vertices of a graph $G$. 
    We say that $w$ resolves $u$ and $v$, if $d(u, w)\neq d(v, w)$. 
    Let $W\subseteq V$ be a set of vertices of $G$ 
     such that every pair of vertices of $G$ is resolved by some vertex in $W$. Then $W$ is called a resolving set of $G$.
      A smallest resolving set of $G$ is called a metric basis, or
      simply a basis, and its cardinality is called the metric dimension of $G$, denoted $\beta(G)$.
\end{definition}
Let $W = \{w_1, w_2\ldots w_n\}$ be a resolving set of a graph $G$. By definition, any distinct vertices $u, v \in G$ are resolved by some vertex $w_i \in W$, meaning they lie in different distance layers of $w_i$. For a vertex $w_i$ and integer $k$, the distance layer $V_{i,k}$ is the set of vertices at distance $k$ from $w_i$. This naturally leads to the notion of the \textit{vector representation} of a vertex $v$ with respect to $W$ defined as $\vr{v} = (d(w_1, v), d(w_2, v),\ldots, d(w_n, v))$. 
It follows that if $W$ is a resolving set, then any two distinct vertices $u, v\in G$ must satisfy $\vr{u}\neq \vr{v}$. 

\bigskip


We now define generalized theta graphs. 
Before proceeding, let $[n]=\{1, 2,\ldots, n\}$.
\begin{definition}
    A generalized theta graph, denoted $\Theta(s_1, s_2, \ldots, s_m)$, 
    is a graph consisting of two endpoints $c_1$ and $c_2$ (referred to as centers) joined by $m$ internally disjoint paths $P_i$ of length $s_i + 1$, where $m,s_i\in \mathbb{N}$ and $i\in [m]$. The integer $m$ is called the multiplicity of the generalized theta graph.

\medskip
    An example of such a graph is given in Fig. \ref{fig:Graph, Generalized Theta Graph}.
\end{definition}
\begin{figure}[h!]
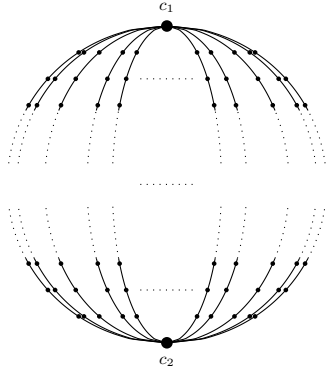

  \centering
  \GraphGeneralizedThetaGraph[0.7]
  \caption{Generalized Theta Graph}
  \label{fig:Graph, Generalized Theta Graph}
\end{figure}

If there are $n_i$ paths of length $s_{i}+1$, then we can denote the generalized theta graph $ \Theta(s_1, \ldots, s_m)$ as 
 $\Theta(s_1^{n_1}, s_2^{n_2}, \ldots, s_k^{n_k})$,
  where  $\sum_{i=1}^{k} n_i = m$. If $s_i=s$ for all $i\in [m]$, then we call such a graph a uniform generalized theta graph and denote it $\Theta(s^m)$. 
 
Let $G = \Theta(s_1, \ldots, s_m)$. For a path $P_i$, in a generalized theta graph, with $i\in [m]$, let $1\leq j \leq s_i$ and let $v_{i,j}$ denote
the vertex of $P_i$, at distance $j$ from the center $c_1$, where the distance is measured along $P_i$. Such a vertex is called an internal vertex. We say that a vertex $v$ is closer to $c_1$ than to $c_2$ if $d(v, c_1) \leq d(v, c_2)$, and closer to $c_2$ than $c_1$ if $d(v, c_2) < d(v, c_1)$. Thus, in the case where $d(v, c_2) = d(v, c_1)$, we assign $v$ to the side of $c_1$ without loss of generality. We will simply write: $v$ is closer to $c_1$ or $v$ is closer to $c_2$, respectively. Furthermore, we define $C_{i,j}$ as the subgraph of $G$ induced by the union of the vertex sets $V(P_i) \cup V(P_j)$.
We refer to $P_i$ as a smallest path if $s_i=s_1$. Finally, let $\tilde{P}_i$ denote the path $P_i$ with both endpoints removed.

\begin{figure}[h!]
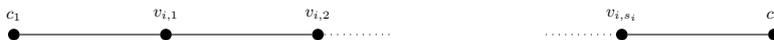

  \centering
  \PathExample[0.6]
  \caption{The vertices on a Generalized Theta Graph within $P_i$}
\end{figure}

Let $v$ be a vertex of $G$. The neighborhood of $v$, denoted $N(v)$, is the set of all
vertices in $G$ adjacent to $v$.
 Given two vertices $u,v\in G$, we say $u$ is
\textit{maximally distant} from $v$, denoted as $u$ MD $v$,
 if for every $w\in N(u)$ it follows that
$d(w, v) \leq d(u, v)$. 
If we have $u$ MD $v$ and $v$ MD $u$, then we say
$u$ and $v$ are \textit{mutually maximally distant}, denoted $u$ MMD $v$. The negation of this statement is denoted as $u$ $\neg$MMD $v$. Lastly, denote by $\mathcal{M}_{G}(v)$ the set of vertices of the graph $G$ that are mutually maximally distant from $v$. Two vertices $u$ and $v$ are called \textit{antipodal} if $d(v,w) = D$, where $D$ is the diameter of the graph $G$. Note that if $u$ and $v$ are antipodal vertices, then $u$ MMD $v$. The following results pertain to cycles. Proofs are omitted when they are straightforward.

\Proposition{ 
\label{prop:MMD pair}
Let $G=C_{2(n + 1)}$ be an even cycle, where $n \geq 1 $. Given a vertex $u\in V(G)$, there exists a unique vertex $v\in V(G)$ such that $ d(u, v) = n + 1$. Furthermore, for each $i \in [n]$, there exists a unique pair of vertices $u', v' \in V(G)$ such that $ d(u, u') = d(u, v') = i$. 
}

\Proposition{ 
\label{prop:WiifTwoMMD}
Let $G=C_n$ and $W\subset V(G)$, where $W=\{w_1, w_2\}$. Then $W$ is not a resolving set if and only if $\mathcal{M}_G(w_1)=\{w_2\}$ and $\mathcal{M}_G(w_2)=\{w_1\}$. 
}
\corollary{
Let $G=C_n$. If $n$ is odd, then any two vertices resolve $G$.
}
\Proposition{
\label{prop:NoTwoDoublePlaces}
Let $G = C_n$ with $n \geq 6$ and let $W = \{w_1, w_2\} \subseteq V(G)$. Then $W$ is not a resolving set if and only if there exist three distinct vertices $u_1, u_2, u_3 \in V(G)$ whose vector representations with respect to $W$ are of the form $(i, i)$, $(j, j)$, and $(k, k)$ for some integers $i$, $j$, and $k$.
}

\begin{proof}
    If $W$ is not a resolving set, then $n$ must be even. Moreover, by Proposition~\ref{prop:MMD pair}, there exist three distinct vertices in $C_n$ whose vector representations are of the form $(i, i)$, $(j, j)$, and $(k, k)$, with the values $i$, $j$, and $k$ not necessarily distinct.
    
    Suppose there exist distinct vertices $u_1, u_2, u_3 \in V(G)$ such that the vector representation of $u_1$ is $(i, i)$, that of $u_2$ is $(j, j)$, and that of $u_3$ is $(k, k)$. If any two of the values in $\{i, j, k\}$ are equal, then $W$ is clearly not a resolving set. Now suppose $i$, $j$, and $k$ are pairwise distinct. Removing the edge between $w_1$ and $w_2$ (in both directions) separates the cycle into two paths, which we denote by $I_1$ and $I_2$. Since there are three vertices but only two paths, at least two of the vertices must lie on the same path. If any two are both in the same path we quickly reach a contradiction, which is always the case, therefore $W$ is not a resolving set.
\end{proof}
\section{Bounds on the Metric Dimension of Generalized Theta Graph}
In this section, we focus on establishing both lower and upper bounds for the metric dimension of generalized theta graphs, as well as presenting a number of complementary results.

\theorem{
\label{thm:UpperBoundGTG}
    Let $G=\Theta(s_1, s_2, \ldots, s_m)$, where $m \geq 3$, then $\beta(G) \leq m$.
}
\begin{proof}
Let $W = \{c_1, w_2, \ldots, w_m\}$ be a set of vertices of the graph $G$, with $w_i=v_{i, \ceil{s_i/2}}$ for $i\in \{2,\ldots,m\}$. We claim that $W$ is a resolving set. Note that for each $k\leq \ceil{s_i/2}$, $|V_{i,k}|=2$. Hence, $V_{i,k}\subseteq V(P_i)$, and therefore internal vertices from different paths have distinct vector representations.

Observe that any two vertices on the path $P_1$ are resolved by the vertex $c_1$. It remains to show that for any path $P_i$, $i\in \{2,\ldots,m\}$, any pair of vertices on $P_i$ are resolved by $W$. We remark that, for each $i$, the distances from $c_1$ and $w_i$ to any vertex $v \in C_{1,i}$ are realized within $C_{1,i}$. Thus it suffices to prove $W_i=\{c_1,w_i\}$ resolves $C_{1,i}$. Because $\nmmd{c_1}{w_i}$, we have that $W_i$ resolves $C_{1,i}$. Repeating this argument for all $i$ completes the proof and shows that $W$ is a resolving set of $G$.
\end{proof}

Having established an upper bound for generalized theta graphs, we will now focus on deriving a lower bound.
\lemma{ 
\label{lem:TwinPathsLemma}
    Let $G$ be a graph in which two distinct vertices $u$ and $v$ are 
    joined by at least two internally disjoint paths, $P_1$ and $P_2$, of equal length, 
    with all internal vertices of these paths having degree 2. Suppose $W\subseteq V(G)$ satisfies $(W \cap V(P_1 \cup P_2)) \setminus \{u, v\} = \emptyset$. Then $W$ is not a resolving set for $G$.
}
\begin{proof}
    Since $P_1$ and $P_2$ are internally disjoint paths of equal length, say $l$, between $u$ and $v$, the graph $G$ is symmetric with respect to these two paths. Thus, the vectors representing any pair of corresponding vertices, such as $v_{1,1}$ and $v_{2,1}$, must be identical. That is, $\vr{v_{1,1}}=\vr{v_{2,1}}$. This implies that $W$ is not a resolving set.
\end{proof}
Two distinct vertices $u$ and $v$ satisfy the $(IP)$ condition if they are connected by a maximal number of internally disjoint paths, say $m \geq 2$, all of which have equal length $s+1$. Furthermore, every internal vertex on these paths must have a degree of 2. We define the set $\mathcal{IP}(G)$ as follows:
\begin{align*}
    \mathcal{IP}(G) = \{ (u, v, s^{m}) \mid u, v \in V(G) \text{ and } u, v \text{ satisfy the } (IP) \text{ condition} \}.
\end{align*}
For uniqueness, we assume that the ordered pair of vertices $(u, v)$ is interchangeable with $(v, u)$. In other words , $(u, v, s^{m})$ and $ (v, u, s^{m})$ refer to the same pair of vertices and paths. We call the set $\mathcal{IP}(G)$ the \textit{Identical Path Set of $G$}. If $|\mathcal{IP}(G)| = n$, we can denote the elements of $\mathcal{IP}(G)$ as $ (u_1, v_1, s_1^{m_1})$, \ldots $ (u_n, v_n, s_n^{m_n})$.
\theorem{[Identical Paths Theorem]
\label{thm:IdenticalPathsTheorem}
Let $G$ be a graph, such that $|\mathcal{IP}(G)| = n$ with $n \geq 1$. Then
\begin{align}
    \beta(G) \geq \sum_{i=1}^n (m_i - 1),
\end{align}
Moreover, if $W \subseteq V(G)$ is a resolving set for $G$, then for each pair of vertices $(u_i, v_i)$ that satisfies the \textbf{($IP$)} condition, $W$ must include at least one internal vertex from each of the $m_i - 1$ paths of length $s_i$ connecting $u_i$ and $v_i$.
}

\begin{proof}
    Let $(u_i, v_i)$ be a pair of vertices in $G$ that satisfies the \textbf{($IP$)} condition, for $1\leq i \leq n$. The paths between $u_i$ and $v_i$ form an induced subgraph of $G$, denoted $G_i$, that is isomorphic to a generalized theta graph, $\Theta(s_i^{m_i})$. In this subgraph,  $u_i$ and $v_i$ act as the centers of a generalized theta graph. Let $W$ be a resolving set of $G$, then $W$ can be expressed as a union of pairwise disjoint subsets, $W=W'\cup W_1 \cup \ldots \cup W_n$, where $W_i$ denotes the subset of resolving vertices from $V(\Theta(s_i^{m_i}))\setminus \{u_i, v_i\}$, 
       while $W'$ denotes the subset of resolving vertices lying outside of these sets, that is $v\in V(G)\setminus \pr{\bigcup_{i=1}^n \pr{V(\Theta(s_i^{m_i}))\setminus \{u_i, v_i\}}}$. For each $W_i$, we claim that $|W_i|\geq m_i-1$.
        Suppose for the sake of contradiction that $|W_i|< m_i - 1$.
        Since $W_i$ corresponds to all resolving vertices in $\Theta(s_i^{m_i})$, this implies that there exists two distinct paths of equal 
        lengths that satisfy Lemma \ref{lem:TwinPathsLemma},
         thus $W_i$ itself cannot be a resolving set. 
         Therefore, it follows that $|W|\geq \sum_{i=1}^n(m_i - 1)$. 
         Furthermore, for at least $m_i-1$ paths in each theta graph, a resolving vertex must be chosen among the internal vertices of each path.
\end{proof}

\corollary{
    \label{cor:LowerBoundUniform}
    Let $G=\Theta(s^m)$, we have that $\beta(G)\geq m-1$.
    Furthermore, if $W$ is a resolving set, then for $m-1$ distinct paths, at least one internal vertex from each path must be included in $W$. 
}

\corollary{
    \label{cor:s1m1s2m2}
    Let $G=\Theta(s_1^{m_1}, s_2^{m_2})$, where $m_1, m_2 \geq 2$. Then $\beta(G)\geq m_1+m_2-2$, furthermore, if $W$ is a resolving set, then for $m_1-1$ distinct paths corresponding to $s_1$ and $m_2-1$ distinct paths corresponding to $s_2$, at least one internal vertex from each path must be included in $W$.
}
\corollary{
\label{cor:biggerThanShortestPaths}
Let $G = \Theta(s_1^{p}, s_{p+1}, \ldots, s_{p+q})$, where $p + q = m$ and $p,q \geq 2$. Then $\beta(G) \geq p+1$.
}
\begin{proof}
    We first show that a resolving set with $|W|=p-1$ does not exist. Fix $p$ and $q$. By Theorem \ref{thm:IdenticalPathsTheorem}, $p - 1$ of the smallest paths must each contain a resolving vertex. Let $W'$ be the set of these resolving vertices. It follows that none of these vertices are mutually maximally distant to any of the centers. There will always remain at least two paths that do not contain any resolving vertices, say $P_i$ and $P_j$. For $u\in V(P_i)\cap N(c)$ and $v\in V(P_j)\cap N(c)$, we have that $r(u\mid W') =r(v\mid W'),$ where $c$ is one of the centers.
    
     Now suppose $|W|=p$, add a new resolving vertex $w$ on any path without a resolving vertex. The issue as before still persists. Therefore $W$ is not a resolving set.
\end{proof}
\corollary{
\label{cor:sallbuts2}
Let $G = \Theta(s_1^{m-1}, s_2)$, where $m \geq 3$. Then $\beta(G) \geq m-1$.
}
\begin{proof}
    The proof follows as the first part of Corollary \ref{cor:biggerThanShortestPaths}.
\end{proof}
Let $G$ be a graph, and let $P = (v_1, v_2, \ldots, v_n)$ be a path in $G$ such that every internal vertex of $P$ (i.e., $v_2, \ldots, v_{n-1}$) has degree two. Let $S = V(P)$, and let $w_1$, $w_2$ be vertices not in $S$. We define $M_i$ to be the set of vertices of $S$ that are mutually maximally distant from $w_i$, for $i = 1, 2$, \text{i.e.,} $M_i := \mathcal{M}_G(w_i) \cap S$.
\proposition{
\label{prop:MMD differen}
 If $M_1 \ne M_2$, then $W = \{w_1, w_2\}$ resolves $P$.
}
\begin{proof}
We prove the contrapositive, if $W = \{w_1, w_2\}$ does not resolve $P$, then $M_1 = M_2$. Suppose $W$ does not resolve $P$. Then there exist two distinct vertices $v_{k_1}, v_{k_2} \in S$ with $k_1 < k_2$ such that $v_{k_1}$ and $v_{k_2}$ have the same vector representation. It follows that
\begin{align}
    d(w_i, v_{k_1 + t}) = d(w_i, v_{k_2 - t})
\end{align}
for all $t$ such that $k_1 + t \leq k_2 - t$ and $i \in \{1, 2\}$. Continuing this process, we eventually reach a midpoint or a pair of vertices $v_r$ and $v_s$ (possibly equal) such that for all $v \in S$,
\begin{align}
    d(w_i, v) < d(w_i, v_r) = d(w_i, v_s),
\end{align}
for each $i \in \{1, 2\}$. Hence, $v_r, v_s \in M_1$ and $v_r, v_s \in M_2$, so $M_1 = M_2$.
\end{proof}

Proposition \ref{prop:MMD differen} guarantees that two vertices in the resolving set uniquely resolve all vertices on a path. It remains to check that vertices from distinct paths, as well as the centers, also admit distinct representations.
\theorem{
    \label{thm:s^m-2,s2,s3}
    Let $G=\Theta(s_1^{m-2}, s_2, s_3)$, where $s_2\leq s_3$ and $m\geq 4$, then $\beta(G)=m-1$.
}
\begin{proof}
    Let $W$ be a resolving set. By  Corollary \ref{cor:biggerThanShortestPaths}, we have $|W|\geq m-1$.

     We claim that $W = \{w_1, \ldots, w_{m-3}, w_{m-1}, w_{m}\}$ is a resolving set, where $w_i = v_{i, \lceil s_1/2 \rceil}$ for $i \in [m - 3]$, and for $j \in \{m - 1, m\}$, we define $w_j = v_{j, \gamma_j}$ with
    \begin{align}
        \gamma_j = \left\lfloor \frac{s_{j - m + 3} + s_1 + 2}{2} \right\rfloor.
    \end{align}

     It follows that $\mmd{w_{m-1}}{c_1}$ in $C_{m-1,1}$ and $\mmd{w_{m}}{c_1}$ in $C_{m,1}$. Then, any distinct pair $i,j\in[m]\setminus\{m-2\}$, it's clear we can view $C_{i,j}$ independently, furthermore $\nmmd{w_i}{w_j}$ in $C_{i,j}$, hence they resolve $C_{i,j}$. 

     Now we show vertices in $\tilde{P}_{m-2}$ are resolved by $W$. Since $V_{i,k}\subseteq V(P_i)$ for $i\in[m-3]$ and $k\in \ceil{s_1/2}$, then the vector representation of vertices in $\tilde{P}_{m-2}$ cannot be the same as those in $V(P_i)$.
     
     Let $u\in V(\tilde{P}_{m-2})$ and $v\in V(\tilde{P}_{m})$, if both vertices have the same vector representations then $d(w_1, u)=d(w_1,v)$. Therefore we consider when both vertices are the same distance to $w_1$. Let $v=v_{m,k}$, if $k>\gamma_{m}$ then it follows that $u$ and $v$ have distinct vector representations since $d(w_m, v) < d(w_m, u)$. If $k<\gamma_m$, then since $d(w_{m-1}, v_{m-2,1}) < d(w_{m-1}, v_{m,1})$ if $u$ and $v$ had the same distance to $w_1$ then it wouldn't be the case that $d(w_{m-1}, u)=d(w_{m-1}, v)$ since going to the left of both paths would never yield vertices with the same vector representation in respect to $w_1$ and $w_{m-1}$. A similar case happens if $v\in V(\tilde{P}_{m-1})$. 
     
     Therefore we only need to show that any pair of vertices in $V(\tilde{{P}}_{m-2})$ have unique vector representation. Since $w_1$ and $w_{m}$ are mutually maximally distant to different vertices in $V(\tilde{{P}}_{m-2})$, by Proposition \ref{prop:MMD differen} each pair has unique vector representations. Hence, $W$ is a resolving set.
 \end{proof}
\theorem{
    \label{thm:onedifferents1>s2}
    Let $G=\Theta(s_1^{m-1}, s_2)$, where $m\geq 4$ and $s_2>s_1$, we have
    \begin{align}
    \beta(G) = \begin{cases}
        m &\, \text{if $G=\Theta(1^{m-1}, 3)$ or $G=\Theta(2^{3}, 4)$,}\\
        m-1 &\, \text{otherwise}
    \end{cases}
    \end{align}
}
\begin{proof}
    By Corollary \ref{cor:sallbuts2}, we have $\beta(G) \geq m-1$. We first address the outlier cases where $\beta(G) = m$. For $G = \Theta(2^3, 4)$, an exhaustive check of all possible resolving sets confirmsthat no set of size $3$ resolves the graph, and thus $\beta(G) = 4$.

Now consider the graph $G = \Theta(1^{m-1}, 3)$. We claim that no resolving set of size $m - 1$ exists. By Theorem \ref{thm:IdenticalPathsTheorem}, and without loss of generality, we may assume the resolving set has the form $W = \{v_{1,1}, \ldots, v_{m-2,1}, w\}$. If $w = v_{m-1,1}$ or $w = v_{m,2}$, then the two centers have the same vector representation. If $w$ is one of the centers, then the neighbors of the center that lie in $P_{m-1}$ and $P_m$ have the same vector representations. If $w = v_{m,1}$, then $\vr{v_{m-1,1}} = \vr{v_{m,3}}$.
The case $w = v_{m,3}$ is similar. Thus, no resolving set of size $m - 1$ exists, and so $\beta(G) \geq m$. By Theorem \ref{thm:UpperBoundGTG}, we have $\beta(G) = m.$

The remainder of the proof is divided into three cases.

\textbf{Case 1}: Let $G=\Theta(2^{m-1},s)$, where $m\geq 5$ and $s>2$. We claim $W = \{w_1, \ldots, w_k, w_{k+1}, \ldots, w_{m-1}\}$ is a resolving set, where $k=\floor{\frac{m-1}{2}}$. If $i\leq k$, then $w_i=v_{i,1}$ otherwise $w_i=v_{i,2}$.

For each $j=1$, we have $V_{i,j}\subseteq V(P_i)$ for $i\in[m-1]$, therefore internal vertices from different paths have distinct vector representations. It also clear that $w_i$ resolves vertices in $\tilde{P}_i$.

By our choice of resolving set, all resolving vertices lie in smallest paths. It follows we can view $C_{i,j}$ independently for distinct $i,j\in[m]$. Suppose both $i,j\leq k$ or both $k<i,j\leq m-1$, then $C_{i,j}$ is resolved by $\{w_i,w_j\}$ since $\nmmd{w_i}{w_j}$ in $C_{i,j}$, therefore the centers and any other vertex have distinct vector representations. This only works because there is always a pair of $i,j$, which is why this resolving set does not work for $G=\Theta(2^3,4)$. 
Lastly, our choice of $w_i$ and $w_j$ are both mutually maximally distant to different vertices in $P_m$, therefore all pairs of vertices in $P_m$ have distinct vector representations by Proposition \ref{prop:MMD differen}. Thus $W$ is a resolving set.

\textbf{Case 2}: Let $G = \Theta(1^{m-1}, s)$. We claim $W=\{v_{1,1}, \ldots, v_{m-2,1}, v_{m,1}\}$ is a resolving set. It is easy to see that $W$ resolves all vertices in $V(G)\setminus V(\tilde{P}_{m})$. It follows that $C_{1,m}$ is resolved by $v_{1,1}$ and $v_{m,1}$ as they are not mutually maximally distant in $C_{1,m}$. Therefore $W$ is a resolving set. This resolving set also works for when $G=\Theta(2^3,s)$ with $s\neq 4$, as one can verify.

\textbf{Case 3}: Let $G = \Theta(s_1^{m-1}, s_2)$ with $s_1 \geq 3$. We claim that $W = \{v_{1,1}, \ldots, v_{m-2,1}, v_{m-1,2}\}$ is a resolving set. Since all resolving vertices lie in the smallest paths, it follows that for any distinct $i, j \in [m]$, we can consider $C_{i,j}$ independently. Since $\nmmd{w_i}{w_j}$ in $C_{i,j}$, both vertices resolve $C_{i,j}$. Therefore, $W$ resolves all vertices in $V(G)\setminus V(\tilde{P}_{m})$.

We now verify that the vertices in $P_m$ are distinguished from all other vertices. Observe that for any vertex in $P_m$, the first $m-2$ coordinates of its vector representations are all equal. However, this is not the case for any vertex in $P_1, \ldots, P_{m-2}$, since each $v_{i,1} \in W$ lies on its respective path. Next, the $(m-1)$th coordinate in the vector representation of $P_{m-1}$ and $P_{m}$ differ when the first $m-2$ coordinates are the same. Finally, to show that the vertices within $P_m$ have distinct vector representations from each other, note that $v_{m-2,1}$ and $v_{m-1,2}$ are each mutually maximally distant to different vertices in $P_m$. Therefore, by Proposition \ref{prop:MMD differen}, all vertices in $P_m$ have unique vector representations. Hence, $W$ is a resolving set.
\end{proof}
\theorem{
\label{thm:Boundfors2<s1OneDifferent}
    Let $G=\Theta(s_1, s_2^{m-1})$, where $m\geq 4$, then
    \begin{align}
        m-2\leq \beta(G) \leq m-1
    \end{align}
}
\begin{proof}
    By Theorem \ref{thm:IdenticalPathsTheorem}, we have that $\beta(G)\geq m-2$,
     thus establishing the lower bound. Let $\gamma = \bfl{\frac{s_1 + s_2 }{2}}+1$,
      we claim that $W=\{c_1, v_{2,\gamma}, \ldots, v_{m-1,\gamma}\}$ is a resolving set.
       Let $w_i= v_{i, \gamma}$, for $2\leq i\leq m-1$. By our choice of $w_i$ it follows that $w_i$ never needs to traverse to vertices in $\tilde{P}_1$ to reach any other vertex. Thus for any two paths $P_i$ and $P_j$, for distinct $i,j\in\{2,\ldots, m-1\}$, we can view $C_{i,j}$ independently. It follows that $\nmmd{w_i}{w_j}$ in $C_{i,j}$ and thus its vertices are resolved. 

       Since $c_1 \in W$, all vertices in $\tilde{P}_1$ have distinct vector representations to each other. Furthermore, for each vertex in $\tilde{P}_1$, the vector representation differs only in the first coordinate, while all other coordinates remain the same. This is not the case for vertices in $\tilde{P}_i$ for $i \in \{2, \ldots, m - 1\}$. As a result, the vertices in each $\tilde{P}_i$ have vector representations that are distinct from those in $V(G) \setminus V(\tilde{P}_m)$.

        Lastly, we show that the vector representations of the vertices in $\tilde{P}_m$ are all distinct. As before, it is clear that these vectors are distinct from those of vertices in $V(G) \setminus V(\tilde{P}_1)$. Additionally, note that $c_1$ and each $w_i$ are mutually maximally distant to different vertices in $V(\tilde{P}_m)$. Therefore, by Proposition \ref{prop:MMD differen}, all vertices in $\tilde{P}_m$ have distinct vector representations from each other. To distinguish vertices in $\tilde{P}_1$ from those in $\tilde{P}_m$, observe that $\nmmd{w_i}{c_1}$ in $C_{i,m}$ for $i \in \{2, \ldots, m-1\}$. Thus, $d(w_i, v_{1,1}) < d(w_i, v_{m,1}).$ It follows that if $u \in V(\tilde{P}_1)$ and $v \in V(\tilde{P}_m)$ had the same vector representation, then their first coordinates would have to agree. However, since the distances from $w_i$ to $u$ and $v$ differ, it cannot be that such $u$ and $v$ exist. Hence, $W$ is a resolving set.
\end{proof}
\theorem{
\label{thm:Distinct,si}
   Let $G = \Theta(s_1, s_2, \ldots, s_m)$, where all $s_i$'s are distinct and $m > 5$. Then $\beta(G) \leq m - 2$.
}

\begin{proof}
     Define the function
    \begin{align}
        \mu(i) = \begin{cases}
        \floor{\frac{s_i+s_1+2}{2}},&\quad \text{if $i$ is even},
        \\
        s_i + 1 -\floor{\frac{s_i+s_1+2}{2}},&\quad \text{if $i$ is odd}.
        \end{cases}
    \end{align}
    We claim $W = \{w_3, w_4,\ldots,w_m\}$ is a resolving set, where $w_i=v_{i,\mu(i)}$. Further define $W_{i,j}=\{w_i, w_j\}$. Essentially, we partitioned the resolving vertices into those that are mutually maximally distant to $c_1$ in $C_{1,i}$, when $i$ is even, and those that are  mutually maximally distant to $c_2$ in $C_{1,i}$, when $i$ is odd. Therefore, we can  consider $C_{i,j}$ independently. 

    If $i$ and $j$ are the same parities, then $\nmmd{w_i}{w_j}$ in $C_{i,j}$ and so $C_{i,j}$ is resolved by $W_{i,j}$. Therefore all vertices in a path have unique vector representations.
    
    Now, let $i$ and $j$ be different parities, say $i$ is even and $j$ is odd. 
    Since all vertices in a path have unique vector representations, it follows that if $u,v\in C_{i,j}$ have the same vector representations then, without loss of generality, $u\in \tilde{P}_i$ and $v\in \tilde{P}_j$. 
    Furthermore, define $I_1=[w_i, c_1, w_j]$ and $I_2=[w_i, c_2,w_j]$. If $u$ and $v$ have the same vector representation, we cannot have both $u,v\in I_1$ or $u,v\in I_2$. We also cannot have $u\in I_2$ and $v\in I_1$, since it will lead to $d(w_i, v) < d(w_i, u)$. Thus $u\in I_1$ and $v\in I_2$.
    Now take $k\in\{3,\ldots,m\}\setminus\{i,j\}$, where $k$ is also even. If $u$ and $v$ have the same vector representations then $d(w_k,u)=d(w_k,v)$. Define 
    \begin{align}
         x_1=d(w_{k}, c_1), \quad
        x_2 = d(w_{k}, c_2),\quad
        y= (c_1, u),\quad z_1=d(c_1,w_j),\quad
        z_2=d(v,c_2).
    \end{align}
    Therefore, we must have $x_1+z_1=x_2+z_1+y+z_2$ and $x_1+y=x_2+z_2$. Combining both equations leads to $y=0$, which means $u=c_1$, a contradiction since for any $\ell\in\{3,\ldots,m\}\setminus\{i,j,k\}$ that is odd we have that $d(w_\ell,c_1)\neq d(w_\ell,v)$.

    For $P_1$ we have that $w_3$ and $w_4$ are mutually maximally distant to different vertices in $P_1$, by Proposition \ref{prop:MMD differen} all pair of vertices in $P_i$ have unique vector representations. The same argument follows for $P_2$. The centers are resolved by the resolving vertices mutually maximally distant to them. By our construction of $W$ and because $s_1$ and $s_2$ are distinct, all vertices in $P_1$ and $P_2$ are resolved by $W$. Hence, $W$ is a resolving set.
\end{proof}
\theorem{
Let $G = \Theta(s_1, s_2, \ldots, s_m)$, such that $s_{i+1} = s_{i}+1$ and $m>6$, then $\beta(G)=m-3$.
}
\begin{proof}
    By Theorem \ref{thm:TotalBound}, the resolving set, $W$, must satisfy $|W| \geq m - 3$. Define the function
    \begin{align}
        \mu(i) = \begin{cases}
        \ceil{\frac{s_i+s_1+2}{2}},&\quad \text{if $i$ is even},
        \\
        s_i + 1 -\ceil{\frac{s_i+s_1+2}{2}},&\quad \text{if $i$ is odd}.
        \end{cases}
    \end{align}
     We claim that $W=\{w_4,w_5,\ldots, w_{m}\}$ is a resolving set, where $w_i=v_{i,\mu(i)}$. The proof follows similar to Theorem \ref{thm:Distinct,si}.
\end{proof}
\lemma{
\label{lem:GeneralizedTwinsPath}
    Let $G$ be a graph and $u,v\in V(G)$, such that 
    \begin{enumerate}
        \item $u$ and $v$ are joined by at least two internally disjoint paths, 
        where the internal vertices have degree 2. Let us call these paths $P_1, P_2$, where the endpoints are both $u$ and $v$.
        \item $d(u, v) +2 < |V(P_1)| ,|V(P_2)| $
    \end{enumerate}
    Let $W\subseteq V(G)$, such that $W\setminus\{u,v\}\cap V(P_1, P_2) = \emptyset$, then $W$ cannot be a resolving set.
}

\begin{proof}
    For sake of a contradiction, suppose $W$ is a resolving set.
     If $|V(P_1)|= |V(P_2)|$ then Lemma $\ref{lem:TwinPathsLemma}$ follows and $W$ cannot be a resolving set. 
     Suppose then that $|V(P_1)| < |V(P_2)|$. Let us use the same 
     notation we have used throughout the paper, namely that 
     $v_{i, j}$ is an internal vertex on the $i$th path and $j$th position from $u$. 
     
     Then we claim that for $w\in W$, 
     we can take $d(w, u)$ and see that $d(w, v_{1, 1})=d(w, v_{2, 1}) =d(w, u) + 1$. 
     To show this first take note that we either have that $v\notin I[w, u]$
      or $v\in I[w, u]$. If we have the first case then our claim is obviously 
      true. Suppose we have the latter then we can split our interval as 
      $I[w, u] = I[w, v]\cup I[v, u]$. We can focus solely on 
      $I[v,u]$, we have that 
      $d(v, u) + d(u, v_{i, 1}) = d(v, u) + 1 \leq d( v, v_{1, 2},v_{1, 1}) < d(v, v_{2, 2},v_{2, 1})$.
      In whatever case $d(w, v_{1, 1}) = d(w, v_{2, 1}) = d(w, v) + d(v, u) + 1$, therefore $W$ is not a resolving set.
\end{proof}
\corollary{
Let $G = \Theta(s_1, s_2, \ldots, s_m)$, where $m > 5$ and $|s_i - s_j| \geq 2$ for all $i \neq j$. Then
\begin{align*}
    \beta(G) = m - 2
\end{align*}
}
\begin{proof}
    Let $W$ be a resolving set, if $|W|< m-2$, then two paths would satisfy the conditions in Lemma \ref{lem:GeneralizedTwinsPath}, it follows $|W|\geq m-2$. By Theorem \ref{thm:Distinct,si}, we have our result.
\end{proof}

\theorem{
\label{thm:LowerBoundGTG}
    Let $G=\Theta(s_1, \ldots, s_m)$, with $m\geq3$, then $\beta(G)\geq m-3$.
}
\begin{proof}
    The cases for when $m\leq 5$ are trivial and so we assume $m>5$. Suppose $W$ is a resolving set such that $|W|<m-3$. Then every resolving set has at least 4 distinct paths such that none of their internal vertices are within $W$. By Lemma \ref{lem:TwinPathsLemma} none of the paths can be of the same length. This leaves us with 4 distinct paths, of which two will always fulfill the conditions of Lemma \ref{lem:GeneralizedTwinsPath}. Hence, $W$ could have not been a resolving set to begin with.
\end{proof}
Therefore, we have established a bound for the metric dimension of theta graphs, which yields Theorem \ref{thm:TotalBound}. We now investigate the metric dimension of uniform theta graphs.

\theorem{
    \label{thm:UGTG1^m}
    Let $G=\Theta(1^m)$, then $\beta(G)=m$.
}
\begin{proof}
    Suppose $W$ is a resolving set, by  Theorem \ref{thm:IdenticalPathsTheorem}, we have $|W|\geq m - 1$. Furthermore, each chosen $m-1$ path contributes a resolving vertex. If $|W|=m - 1$,
    without loss of generality, let $W=\{v_{1, 1}, \ldots, v_{m-1, 1}\}$. Then $\vr{c_1}=\vr{c_2}=(1,1, \ldots, 1)$. Hence, $W$ is not a resolving set, by Theorem \ref{thm:UpperBoundGTG}, $\beta(G)=m$.
\end{proof}
\theorem{
    \label{thm:UGTG2^m}
    Let $G=\Theta(2^m)$, with $m\geq 3$, we have the following
    \begin{align}
        \beta(G) = \begin{cases}
            m &\, \text{ if $m\leq 4$}\\
            m - 1 &\, \text{ otherwise}
        \end{cases}
    \end{align}
}
\begin{proof}
    By Theorem \ref{thm:IdenticalPathsTheorem}, if $W$ is a resolving set, 
    then $|W|\geq m - 1$ and each path must contribute a resolving vertex. 
     Suppose $|W|=m-1$, with $W = \{v_{1, k},v_{2, k},\ldots, v_{m-1, k} \}$, without loss of generality suppose $k=1$, we immediately see that $\vr{v_{m, 1}}=\vr{c_2}=(2, 2,\ldots, 2)$. Thus $W$ is not a resolving set.
    
    Let $m = 3$ and $W =\{v_{i, 1}, v_{j, 2}\}$. Without loss of generality, $W = \{v_{1, 1}, v_{2, 2}\}$. Notice however that $\vr{c_1}=\vr{v_{1,2}}=(1,2)$.
    If $m = 4$, let $W =\{v_{i, 1}, v_{j, 2}, v_{k, \ell}\}$. Without loss of generality let $W = \{v_{1, 1}, v_{2, 2}, v_{3, 2}\}$. Notice however that $\vr{c_1}=\vr{v_{1,2}}=(1,2, 2)$. In both cases we have that $W$ is not a resolving set. Then by Theorem \ref{thm:UpperBoundGTG}, $\beta(G)=m$, for $m\in\{3,4\}$.
    
    Let $m>4$. We claim $W=\{w_1,\ldots,w_{m-1}\}$ is a resolving set, where $w_{i}=v_{i,1}$ if $i\leq \floor{\frac{m+1}{2}}$, otherwise $w_i=v_{i,2}$. Since $V_{i,k}\subseteq V(P_i)$ for $i\in[m-1]$ and $k=1$, it follows any two vertices in different paths have unique vector representations. It is easy to see that the rest of the vertices have unique vector representations. Therefore $W$ is a resolving set.
\end{proof}
\theorem{
\label{thm:Uniform Theta Proof}
    Let $G=\Theta(s^m)$, where $m \geq 3$ and $s \geq 1$, we have 
    \begin{align}
        \beta(G)=\begin{cases}
            m-1 &\, \text{if $m\geq 5$ and $s\geq2$}\\
            m-1 &\, \text{if $m=4$ and $s>2$}\\
            m &\, \text{otherwise}
        \end{cases}
    \end{align}
}
\begin{proof}
    Suppose $m\geq 4$ and $s>2$. We claim $W=\{w_1,\ldots,w_{m-1}\}$ is a resolving set, where $w_i=v_{i,1}$ for $i\in[m-2]$ and $w_{m-1}=v_{m-1,2}$. Since all paths are the same lengths, for distinct $i,j\in[m-1]$, we may view $C_{i,j}$ independently. Clearly, $\nmmd{w_i}{w_j}$, therefore $C_{i,j}$ are resolved by both vertices. 

    We now just need to show that all vertices in $P_m$ have unique vector representations. Since $w_1$ and $w_{m-1}$ are mutually maximally distant to different vertices in $P_m$, by Proposition \ref{prop:MMD differen}, all pairwise vertices in $P_m$ have unique vector representations. Next, notice the first $m-1$ coordinates in the vector representation of $u\in P_m$ are the same, a property that does not hold in $v\in P_i$, where $i\in[m-2]$. 
    
    We now only have to show vertices in $P_{m-1}$ have unique vector representations to those in $P_m$. Let $u\in P_{m-1}$ and $v\in P_{m}$, if both vertices have the same vector representation, then $d(w_1, u)=d(w_1,v)$, meaning $u=v_{m-1,j}$ and $v=v_{m,j}$ for some $j\in[s]$. However, it is always the case that $d(w_{m-1}, u)<d(w_{m-1},v)$. Therefore $u$ and $v$ have distinct vector representations. Hence, $W$ is a resolving set.
     All other cases follow by Theorems \ref{thm:UGTG1^m}, \ref{thm:UGTG2^m} and \ref{thm:GTGEndResult}.
\end{proof}

\section{Metric Dimension of Generalized Theta Graph of multiplicity 3}
In this section, we determine the metric dimension of all theta graphs, namely the generalized theta graphs of multiplicity 3, and show that it is always 2 except for two cases, thus yielding another proof in the results previously obtained in \cite{Knor2022remarks}.
\theorem{
\label{thm:GTG3DiffParitity}
    Let $G = \Theta(s_1, s_2, s_3)$ be a theta graph that is neither of the form $\Theta(s^3)$ nor $\Theta(s_1^2, s_2)$ with $s_2 = s_1 + 2$. Suppose $s_1$, $s_2$, and $s_3$ do not all share the same parity. Then $\beta(G) = 2$.
}

\begin{proof}
    Without loss of generality, assume $s_1$ is odd and $s_2$ is even. We claim that $W = \{v_{1, \bc{ s_1 / 2}}, v_{2, s_2 / 2} \}$ is a resolving set. By our choice of resolving vertices we can view $C_{1,2}$ independently, thus $W$ resolves all the vertices in $C_{1,2}$. Since $s_1$ and $s_2$ are of different parities we have that they cannot be mutually maximally distant to the same set of vertices in $P_3$ and thus all pair of vertices in $P_3$ have unique vector representations by Proposition \ref{prop:MMD differen}. Additionally, $V_{i, j}\subseteq V(P_i)$ for 
    $j\in \{1,\ldots, \bc{s_i/2}\}$ and $i\in [2]$, thus all vertices in $P_1$ and $P_2$ have distinct vector representations to any $v\in V(P_3)\setminus \{c_1, c_2\}$. Thus $W$ is a resolving set.
\end{proof}
\theorem{
 \label{thm:GTG3SameParitity}
 Let $G = \Theta(s_1, s_2, s_3)$ be a theta graph that is neither of the form $\Theta(s^3)$ nor $\Theta(s_1^2, s_2)$ with $s_2 = s_1 + 2$. Suppose $s_1$, $s_2$, and $s_3$ all share the same parity. Then $\beta(G) = 2$.
}
\begin{proof}
We will consider two separate cases.

Suppose $s_1 < s_2 \leq s_3$, and define $\gamma_i = \frac{s_1 + s_i + 2}{2}$. Let $W = \{w_1, w_2\}$, where $w_i = v_{i+1, \gamma_{i+1}}$. We claim that $W$ is a resolving set. First, consider $C_{1,2}$, by construction, $\mmd{w_1}{c_1}$, and $w_2$ forces $c_2$ to act as a resolving vertex within $C_{1,2}$. Since $\nmmd{w_1}{w_2}$ it follows that every pair of vertices in $C_{1,2}$ has a unique vector representation. A similar argument applies to $C_{1,3}$. Finally, observe that in $C_{2,3}$ each resolving vertex does not need to traverse through $P_1$ to reach a vertex in $C_{2,3}$, therefore we may consider it independently. Since both resolving vertices are not mutually maximally distant, all pairs of vertices in $C_{2,3}$ have distinct vector representations. Hence, $W$ is a resolving set.

Suppose $s_1=s_2$ and $s_1+2<s_3$. Then we claim $W=\{v_{2,1},v_{3, 1}\}$ is a resolving set, where $w_1= v_{2,1}$ and $w_2=v_{3, 1}$. One can see that $w_i$ does not need to traverse through $P_1$ to reach a vertex in $P_3$, then we can view $C_{2,3}$ independently. Since $\nmmd{w_1}{w_2}$ all vertices in $C_{2,3}$ have unique vector representations. All vertices in $V(\tilde{P}_1)$ also have vector representations where both coordinates are the same, which does not happen for any other vertices except $c_1$, however one can see that they are unique. Therefore $W$ is a resolving set.
\end{proof}
\theorem{
    \label{thm:GTG3Exceptions}
 Let $G = \Theta(s_1, s_2, s_3)$ be a theta graph that is in the form $\Theta(s^3)$ or $\Theta(s_1^2, s_2)$ with $s_2 = s_1 + 2$. Then $\beta(G) = 3$.
}
\begin{proof}
    We will show that neither of these graphs have a resolving set of size 2. 

    Suppose $G=\Theta(s^3)$. By Lemma \ref{lem:TwinPathsLemma}, the resolving vertices must come from distinct paths. Without loss of generality, let $W=\{v_{1, k}, v_{2, \ell}\}$. If $k=\ell$, both resolving vertices are mutually maximally distant to the same vertex, $u\in P_3$, and thus the neighbors of $u$ must have the same vector representation.
     Then $k\neq \ell$, we can further assume $k<\ell$. If $k+\ell = s+1$, then this means $\mmd{w_1}{w_2}$, therefore $C_{1,2}$ is not resolved by $W$. If $k+\ell < s+1$, then consider the vertices $u=v_{1,s+1+k-\ell}$ and $v=v_{3, s+1-(k+\ell)}$, then 
    \begin{align}
             \vr{u} = \vr{v} = (s+1-\ell, s+1-k).
    \end{align}
    If $k+\ell> s+1$, consider the vertices $u=v_{2,\ell-k}$ and $v=v_{3, 2s+2 - (\ell+k)}$, then
    \begin{align}
             \vr{u} = \vr{v} = (\ell, k).
    \end{align}
    In every case, $W$ cannot be a resolving set.

    Suppose $G=\Theta(s_1^2, s_2)$, with $s_2=s_1+2$. Let $W=\{w_1,w_2\}$, by Lemma \ref{lem:TwinPathsLemma}, one of the resolving vertices must come from $P_1$ or $P_2$. Without loss of generality, let $w_1\in P_1$. If $w_2\in V(P_1)$ or is a center, then it is always the case that there exists neighbors of a center that have the same vector representations.

    Therefore $w_2$ cannot be a center or in $P_1$. Suppose $w_2\in P_2$, and let $W=\{v_{1, k}, v_{2, \ell}\}$ with $w_1=v_{1, k}$ and $w_2= v_{2, \ell}$. If $k=\ell$, we have that both resolving vertices are antipodal to some $u\in P_3$, and so the neighbors of $u$ must have the same vector representation. Without loss of generality, let $k<\ell$. If $k+\ell=s_1+1$, it follows that $\mmd{w_1}{w_2}$ therefore $C_{1,2}$ is not resolved by $W$. If $k + \ell< s_1 + 1$, consider $v=v_{1, s_1 +1 +k - \ell}$ and $u=v_{3, s_1 + 1 - (k + \ell)}$, then
    \begin{align}
        \vr{u}=\vr{v} = (s_1+1-\ell, s_1+1-k).
    \end{align}
     If $k + \ell > s_1 + 1$, consider $v=v_{2, \ell-k}$ and $u=v_{3, 2s_1 + 4 - (k + \ell)}$. Given that $k>1$, we have the
    \begin{align}
        \vr{u}=\vr{v} = (\ell,k).
    \end{align}
    Therefore, if $w_2\in V(P_2)$, then $W$ is not a resolving set.
    
    Lastly, suppose $w_2\in V(P_3)$, then let $W=\{v_{1,k},v_{3, l}\}$. We consider $k = \ell$. If $k\leq \bc{\frac{s_1}{2}}$, consider $u=v_{2, s_1-2k+2}$ and $v=v_{3, s_1 + 2}$, we have
    \begin{align}
        \vr{u}=\vr{v} = (s_1 - k + 2, s_1 - k + 2).
    \end{align}
    If $k=\ceil{s_2/2}$ and $s_2$ is even, it follows that $\mmd{w_1}{w_2}$ therefore $W$ cannot be a resolving vertex. If $k>\bc{\frac{s_1}{2}}$, consider $u=v_{1,1}$ and $v=v_{2, 2s_1 -2k+3}$, we have 
    \begin{align}
        \vr{u}=\vr{v} = (k-1,k+1).
    \end{align}
    It follows that $W$ is not a resolving set. 
     Now we consider the case $k< \ell$.
     \begin{itemize}
         \item[] \textbf{Case 2.3.1} Suppose $k+\ell < s_1 + 2$. Observing $C_{1,2}$, by Proposition \ref{prop:MMD pair} we have that there is two distinct vertices $u$ and $v$ in $C_{1,2}$ such that $d(w_1, u)=d(w_2, v) = s_1 + 2 -j$. By our construction, we always have $u\in V(P_2)\setminus \{c_1, c_2\}$.
         
         If $v=c_2$, then $d(c_1, w_1) + d(w_1, c_2) = s_1+1\implies k+1 = \ell$. Furthermore this means $u = v_{2, s_1+1-2k}$. One can check that $d(w_2, u)=d(w_2, v)= s_1 -k+2$.
         If $v\in V(P_1)\setminus\{c_2\}$, then $v=v_{1, s_2 - \ell+k}$ and $u = v_{2, s_1 + 2 - (\ell+k)}$. One can check that $d(w_2, u)=d(w_2, v)= s_1 -k+2$.
        
        In both cases we have $\vr{u}=\vr{v}$.
        \item[] \textbf{Case 2.3.2} Suppose $k+\ell > s_1 + 2$. If $k + 1 = \ell $, then consider $u = c_1$ and $v = v_{2, 2s_1 + 3 - (l+k)}$. We have the following distances
        \begin{align}
            \vr{u} &= (k, l)\\
            d(w_1, v) &= \min\{d(w_1, c_1, v), d(w_1, c_2, v)\} \\
            &= \min\{2s_1 + 3 - \ell,  \ell-1 \} = \ell-1 = k\\
            d(w_2, v) &= \min\{d(w_2, c_1, v), d(w_2, c_2, v)\} \\
            &= \min\{2s_1 + 3 - k,  k +1 \} = k+1 = \ell
        \end{align}
        And so $\vr{u}=\vr{v}$. If $k+1\neq \ell$, then consider $u=v_{2, 2s_1+3 - (\ell + k)}$ and $v = v_{3, \ell - (k + 1)}$. We have the following distances
        \begin{align}
            d(w_1, u) &= \min\{d(w_1, c_1, u), d(w_1, c_2, u)\} \\
            &= \min\{2s_1 + 3 - \ell,  \ell -1 \} = \ell-1 \\
            d(w_2, u) &= \min\{d(w_2, c_1, u), d(w_2, c_2, u)\} \\
            &= \min\{2s_1 + 3 - k,  k +1 \} = k+1\\
            d(w_1, v) &= \min\{d(w_1, c_1, v), d(w_1, c_2, v)\} \\
            &= \min\{\ell-1, 2s_1+5 - l\} = \ell-1 \\
            d(w_2, v) &= \ell - (\ell - (k + 1)) = k + 1
        \end{align}
        We again get $\vr{u}=\vr{v}$.
        \item[] \textbf{Case 2.3.3} Suppose $k + \ell = s_1 + 2$. If $s_1$ is odd, and $k = \bc{\frac{s_1}{2}}$ and $\ell = \bc{\frac{s_2}{2}}$, then clearly $\vr{c_1}=\vr{c_2}$. Otherwise, consider $u = c_2$ and $v=v_{3, \ell -(k+1)}$. We have the following distances
        \begin{align}
            d(w_1, u) &= (s_1 + 1) - k = (k+\ell - 1) + k = \ell - 1 \\
            d(w_1, v) &= \min\{d(w_1, c_1, v), d(w_1, c_2, v)\}\\
            &= \min\{\ell-1, \ell + 2k + 1\} = \ell-1
            \\
            d(w_2, u) &= (s_1+3) - \ell = k + \ell + 1 - \ell = k + 1 
            \\
            d(w_2, v) &= \ell - (\ell - (k+1)) = k + 1
        \end{align}
        And so $\vr{u}=\vr{v}$.
     \end{itemize}
     Thus $W$ is not a resolving set when $|W|=2$.
\end{proof}
\theorem{
    \label{thm:GTGEndResult}
    \begin{align}
        \beta(\Theta(s_1, s_2, s_3)) = \begin{cases}
            3 &\, \text{if (1) $s_1=s_2=s_3=s$ or (2) $s_1=s_2$ and $s_3 = s_1 + 2$}\\
            2 &\, \text{otherwise}
        \end{cases}
    \end{align}
}
\begin{proof}
    In the first case, Theorem \ref{thm:GTG3Exceptions} establishes that the metric dimension is not 2, combined with Theorem \ref{thm:TotalBound}, this gives the desired result. The next case follows directly from Theorems \ref{thm:GTG3DiffParitity}, and  \ref{thm:GTG3SameParitity}.
\end{proof}
\section{Metric Dimension of Generalized Theta Graphs of multiplicity 4}
In this section, we investigate the metric dimension of the generalized theta graphs of multiplicity 4 through a case-by-case study..

\theorem{ \cite{Khuller1996LandmarksIG}
\label{thm:degreelessthan3}
    Let $G = (V, E)$ be a graph with metric dimension $2$ and let $\{u, v\} \subset V$ be a metric basis of $G$. The following are true:
\begin{itemize}
  \item[1.] There is a unique shortest path $P$ between $u$ and $v$.
  \item[2.] The degrees of $u$ and $v$ are at most $3$.
  \item[3.] Every vertex of $P$, other that $u$ and $v$, has degree at most $5$.
\end{itemize}
}
\corollary{
\label{cor:4NoCenters}
 Let $G=\Theta(s_1,s_2,s_3,s_4)$. No center of $G$ belongs to a resolving set of size 2.
}
Not that the above result also holds for $G=\Theta(s_1,s_2,s_3,s_4, s_5)$

\proposition{
\label{prop:NotMiddle4}
 Let $G=\Theta(s_1,s_2,s_3,s_4)$. If $W$ is a resolving set of $G$ with $|W| = 2$, then the two vertices in $W$ must be closer in distance to the two centers; that is, one vertex must be closer to $c_1$ and the other to $c_2$.
}
\begin{proof}
Let $W = \{w_1, w_2\}$. Suppose, without loss of generality, that both $w_1$ and $w_2$ are closer to $c_1$. Whether the vertices $w_1$ and $w_2$ lie on the same path of $G$ or on different paths, there exist two vertices on distinct paths of $G$, disjoint from those containing $w_1$ and $w_2$, that share the same vertex representation with respect to $W$. Hence, $W$ cannot be a resolving set of $G$.
\end{proof}

\lemma{
\label{lem:Shortest2for4}
    Let $G=\Theta(s_1,s_2,s_3,s_4)$ and let $W \subset G$ be a resolving set such that $|W|=2$.
    If $P$ is a shortest path, it follows that $W\cap V(P) = \emptyset$.
}
\begin{proof}
    Assume $W \cap V(P) \neq \emptyset$. By Corollary~\ref{cor:4NoCenters}, $W$ does not contain any centers. We may write $W = \{w_1, w_2\}$, where $w_1 \in P$ and $w_2$ is an internal vertex.

    Without loss of generality, suppose $d(w_2, c_1) \leq d(w_2, c_2)$, so $w_2$ is closer to $c_1$. Since $w_1 \in P$, and every vertex on $P$ lies on a shortest path between $c_1$ and $c_2$, $w_1$ cannot be mutually maximally distant with either center. Moreover, since $w_2$ is closer to $c_1$, $w_2$ cannot be mutually maximally distant with $c_1$. It follows that there exists a pair of vertices in the neighborhood of $c_1$ whose vector representations are the same. Therefore, $W$ is not a resolving set.
\end{proof}

\theorem{
    Let $G=\Theta(s_1,s_2,s_3,s_4)$, where $s_1=s_2=s_3 <s_4$. Then
    \begin{align*}
    \beta(G) = \begin{cases}
        4 &\, \text{for $\Theta(1^3, 3)$ or $\Theta(2^3, 4)$}\\
        3 &\, \text{otherwise}
    \end{cases}
    \end{align*}
}
\begin{proof}
    This result is a special case of Theorem \ref{thm:onedifferents1>s2}.
\end{proof}
\theorem{
\label{thm:MD43same1dif}
    Let $G = \Theta(s_1, s_2, s_3, s_4)$, such that $s_1 <s_2=s_3=s_4$. Then $\beta(G) = 3.$
}
\begin{proof}
   By Theorem \ref{thm:Boundfors2<s1OneDifferent}, it suffices to show that the graph $G$ does not admit a resolving set $W$ of size 2.
   
By Lemma \ref{lem:Shortest2for4}, $W$ does not contain vertices in $V(P_1)$. Furthermore, the Identical Paths Theorem implies that the two vertices in $W$ must lie on distinct paths of $G$. 

Without loss of generality, let $W = \{w_1, w_2\} = \{v_{3,i}, v_{4,j}\}$. By Proposition~\ref{prop:NotMiddle4}, $w_1$ and $w_2$ cannot both be closer to the same center. Therefore, we can assume $i<\ceil{s_3/2} < j$, making $w_1$ closer to $c_1$ and $w_2$ closer to $c_2$. Consider the sets $I=\mathcal{M}_{C_{1,3}}(w_1)$ and $J=\mathcal{M}_{C_{1,4}}(w_2)$.

It follows that $I$ can contain either only vertices from $P_1$, only vertices from $P_3$, the single vertex $c_2$, the pair $\{v_{1, s_1}, c_2\}$, or the pair $\{v_{3, s_3}, c_2\}$. Similarly, $J$ can contain either only vertices from $P_1$, only vertices from $P_4$, the single vertex $c_1$, the pair $\{v_{1, s_1}, c_1\}$, or the pair $\{v_{4, 1}, c_1\}$.

It suffices to consider $I$, as the case for $J$ is analogous. If $I$ only contains vertices in $P_1$, then $v_{1, s_1}$ and $v_{2, s_2}$ have the same vector representations. If $I$ only contains vertices in $P_3$, then $v_{2, s_2}$ and $v_{3, s_3}$ have the same vector representations.  If $I=\{c_1\}$, then clearly $v_{1,s_1}$ and $v_{3, s_3}$ have the same vector representation. 

In the remaining two cases involving the set $I$, the contradiction cannot be obtained without considering the set $J$ as well. If  $J$ is not in its last two cases, then clearly $J$ produces two vertices with the same vector representations. Therefore, suppose $J$ must be in its last two cases. Moreover, if these last two cases arise, then it follows that $s_1$ and $s_3$ are of different parity. 

If both $I$ and $J$ are in their third option, then the shortest path plays no role in determining distances within the induced subgraph on the vertex set $V(P_2) \cup V(P_3) \cup V(P_4)$ are determined in the centers and rest of the graph. Hence, by Theorem \ref{thm:GTG3Exceptions}, the two vertices cannot resolve this induced subgraph.

If both $I$ and $J$ are in their fourth option, then the distance between $w_1$ and each vertex of $I$ is $(s_3+s_1+1)/2$; similarly for $w_2$ and $J$. Hence, it follows that 
\begin{align}
    i= \frac{s_3+s_1+1}{2}-(s_1+1),
\quad j =  \frac{s_3+s_1+1}{2}+1,
\end{align}
since $i<s_3-s_1$.
Consider the vertices $u=v_{3, s_3-s_1}$ and $v=v_{4, 1}$. One can compute to verify 
\begin{align}
    d(w_1, u) &= s_3 - s_1 - i = i + 1 = d(w_1, v), \\
    d(w_2, u) &= s_3 - j + 1 + s_3 - (s_3 - s_1) + 1 = j - 1 = d(w_2, v).
\end{align}
Therefore both vertices have the same vector representations.

Lastly, the final two cases lead to analogous arguments. Thus, we may take $I$ to be in its third option and $J$ in its fourth option. As in the previous case, the distances from the resolving vertices to the vertices in their respective sets ($I$ and $J$) can be determines, and we obtain 
\begin{align}
    i= \frac{s_3+s_1+1}{2}-s_1,
\quad j =  \frac{s_3+s_1+1}{2}+1.
\end{align}
Consider the vertices $u=v_{2,s_3}$ and $v_{4,s_1+1}$. A simple calculation shows 
\begin{align}
    d(w_1, u) &= s_3-i+2= i + s_1+1 = d(w_1, v), \\
    d(w_2, u) &= j - s_1 - 1= s_3-j + 2 = d(w_2, v).
\end{align}
Therefore both vertices have the same vector representations.

Hence, $W$ is not a resolving set in any case. Hence, and so $|W| \neq 2$. The theorem follows.
\end{proof}
\theorem{
    Let $G = \Theta(s_1^2, s_2, s_3)$, such that $s_2 \leq s_3$, then $\beta(G) = 3$.
}
\begin{proof}
    This is a special case of Theorem \ref{thm:s^m-2,s2,s3}.
\end{proof}
\theorem{
\label{thm:4,s1+1<s2sup}
    Let $G=\Theta(s_1, s_2, s_3, s_4)$, where $s_1+1<s_2$. Then $\beta(G)=3$.
}
\begin{proof}
     Let  $W=\{w_1,w_2\}$ be a resolving set. By Corollary \ref{cor:4NoCenters}, $W$ cannot include any center. By Lemma \ref{lem:Shortest2for4} the shortest paths must not contain resolving vertices. 
     By Theorem \ref{thm:MD43same1dif} we do not consider the case where $s_2=s_3=s_4$. 
    By Lemma \ref{lem:GeneralizedTwinsPath} we cannot have $W\subseteq V(P)$ for any path $P$. Thus, it remains to consider the case when resolving vertices lie on different paths.
     By Proposition \ref{prop:NotMiddle4}, each vertex must be closer to a different center.
     
    For each $k\in \{1,2\}$, we write $w_{k}=v_{i_k, j_k}$, where $i_k\in\{2,3,4\}$ and $j_k\in[s_{i_k}]$. Further, suppose that $w_k$ is closer to $c_k$. By Theorem \ref{thm:degreelessthan3} there is a unique shortest path from $w_1$ to $w_2$, therefore if $s_{i_1}$ and $s_{i_2}$ have the same parity, then it is not possible to have $w_1$ MMD $w_2$.

    There are three cases for the shortest unique path between $w_1$ and $w_2$: $I[w_1, c_1,v_{i_2,1},w_2]$, $I[w_1, c_1, v_{1,1},w_2]$, and $I[w_1, v_{i_1,s_{i_1}},c_2,w_2]$. Since the third case follows by the same reasoning as the first, it is enough to analyze the first two cases. 

    As a first case, suppose the unique shortest path is $I[w_1, c_1,v_{i_2,1},w_2]$. Let $J_k=\mathcal{M}_{C_{1,i_k}}(w_k)$, for $k\in\{1,2\}$. If $J_2\subseteq V(P_{i_2})$, then the unique shortest path cannot be $I[w_1, c_1,v_{i_2,1},w_2]$ as 
    \begin{align}
        d(w_2,v_{i_2,j_2-1},c_1)>d(w_2,v).
    \end{align}
    If $J_2\subseteq V(\tilde{P}_1)$, then there exists a distinct pair of vertices in $N(c_1)$ that have the same vector representation. Finally, suppose $J_2=\{c_1,v_{1,1}\}$. In this case we also examine $J_1$. If $J_1=\{c_2,v_{1,s_1}\}$, then $s_{i_1}$ and $s_{i_2}$ have the same parity since this would mean that both $C_{1,s_{i_1}}$ and $C_{1,s_{i_2}}$ are odd cycles as $w_1$ and $w_2$ are both MMD to two vertices each. Consider \figref{T7}.
    \begin{figure}[hbt!]
  \centering
  \begin{tikzpicture}[scale=0.3, transform shape]  
    \input{T4}             
\filldraw[fill=white, draw=black, line width=0.8pt] (v4_2) circle (6pt);
\node[yshift=18pt,  xshift=5pt, scale=2.5] at (v4_2) {$w_1$};

\filldraw[fill=white, draw=black, line width=0.8pt] (v2_2) circle (6pt);
\node[yshift=-18pt,  xshift=10pt, scale=2.5] at (v2_2) {$w_2$};

\filldraw[fill=black, draw=black, line width=0.8pt] (v1_1) circle (6pt);
\node[yshift=18pt,  xshift=-10pt, scale=2.5] at (v1_1) {$u$};

\filldraw[fill=black, draw=black, line width=0.8pt] (v1_9) circle (6pt);
\node[yshift=18pt,  xshift=10pt, scale=2.5] at (v1_9) {$v$};
  \end{tikzpicture}
  \caption{\label{fig:T7}}
\end{figure}

    Furthermore it follows that $P_1$ does not need to be traversed and that $w_1$ is MMD to one such vertex $v\in V(P_{s_{i_2}})$, such that $d(w_1,v)=\frac{s_3+s_4}{2}+1$. Notice that $d(w_1,c_1)=\frac{s_3-s_1-1}{2}+1$ and $d(w_2,c_1)=\frac{s_4+s_1+1}{2}$, thus $d(w_1,c_1,w_2)= \frac{s_3+s_4}{2}+1$. Therefore $v=w_2$ and it follows that the shortest path is not unique. The remaining cases follow analogously by considering $J_1$ instead of $J_2$. 

     Now consider the case where the unique shortest path is $I[w_1, c_1, v_{1,1},w_2]$. If $J_2\subseteq V(P_1)$, then the shortest path cannot be $d(w_1,c_1,v_{1,1},w_2)$ as $d(w_1,c_1,v_{i_2,1},w_2)$ is shorter. If $J_2\subseteq V(\tilde{P}_{i_2})$, then there exists a pair of vertices in $N(c_1)$ having the same vector representation.
     
     Finally, consider the case $J_2=\{c_1, v_{i_2,1}\}$. In this situation we also examine $J_1$. If $J_{1}=\{c_2,v_{i_1,s_{i_1}}\}$, then by the same reasoning as before $s_{i_1}$ and $s_{i_2}$ have the same parity, however we may not consider $C_{i_1,i_2}$ independently as traversing through $P_1$ is necessary. Nevertheless, both $w_1$ and $w_2$ are still determined. Let $r_k = (s_{i_k}+s_1 + 1)/2$, then $w_1= v_{i_1, s_{i_1}-r_1}$ and $w_2= v_{i_2, r_2+1}$. Consider the vertices $u=v_{i_1, s_{i_1}}$ and $v=v_{i_2,2r_1-s_{i_1}}$. Then, we obtain
    \begin{align}
        \vr{u} = (r_1, r_2-s_1) = \vr{v}.
    \end{align}
    The remaining cases follow by considering $J_1$ instead of $J_2$. Therefore $W$ is not a resolving set.

        Thus, there does not exist a resolving set of size 2. We have that $W= \{c_1, v_{3, \gamma_3}, v_{4, \gamma_4}\}$ is a resolving set, where $\gamma_i=\floor{\frac{s_1+s_{i}+2}{2}}$. Since $s_2> s_1+1$, this ensures that the vector representations of the vertices in $P_1$ and $P_2$ are distinct.
     The remainder of the proof proceeds as in Theorem \ref{thm:s^m-2,s2,s3}. 
\end{proof}
\theorem{
    
    Let $G=\Theta(s_1, s_2, s_3, s_4)$, where $s_1+1= s_2=s_3$ and $s_4-s_1 \in \{2,3\}$. Then $\beta(G)=3.$
}
\begin{proof}
    We start by showing that no resolving set $W=\{w_1,w_2\}$ exists. If it did, by Corollary \ref{cor:4NoCenters} the centers are not in $W$. By Proposition \ref{prop:NotMiddle4}, let's suppose $w_i$ is closer to $c_i$. By Lemma \ref{lem:Shortest2for4}, no resolving vertex is from $V(P_1)$. Since $s_2=s_3=s_1+1$, we have that $w_i$ never needs to traverse through $V(P_1)$ to reach a vertex $V(G)\setminus V(P_1)$. Therefore, if $s_4=s_1+3$ then by Theorem \ref{thm:GTGEndResult}, $W$ never resolves $G$.  Then we only consider $s_4=s_1+2$. 
    By Lemma \ref{lem:TwinPathsLemma}, a resolving vertex must lie on $P_2$ or $P_3$.

    Consider the case when both resolving vertices are in $P_2$ or $P_3$, without loss of generality let them be in $P_3$. By our assertion that $w_i$ does not need $P_1$, we get that $P_3$ ``acts" like a smallest path. Therefore $\vr{v_{2,1}} =\vr{v_{4,1}}$.

    Now suppose $w_2\in P_2$ instead. If $w_1$ MMD $w_2$ we have that $C_{2,3}$ is not resolved as it is an even cycle. If $w_1$ $\neg$MMD $w_2$, then both resolving vertices are not neighbors of the centers they are close to, thus $\vr{v_{1,1}}= \vr{v_{4,1}}$ or $\vr{v_{1,s_1}}= \vr{v_{4,s_4}}$.

     Since $s_1$ and $s_4$ have the same parity, it must be $w_2$ is MMD to $c_1$ or a vertex in $P_1$ within $C_{1,4}$. If $w_2$ is MMD to $c_1$, since $w_1$ is closer to $c_1$, we have $\vr{v_{1,1}}=\vr{v_{4,1}}$. If $w_2$ is MMD to a vertex in $P_1$, since $w_1$ is closer to $c_1$, we have $\vr{v_{1,1}}=\vr{v_{2,1}}$.
    Thus no resolving set of size 2 exists.

      We have that $W= \{c_1, v_{3, \gamma_3}, v_{4, \gamma_4}\}$ is a resolving set, where $\gamma_i=\floor{\frac{s_i+s_{1}+2}{2}}$. By our construction, it follows $\gamma_3=s_3$, therefore we have $v_{3,\gamma_3}=w_2$ is MMD to $v_{2,1}$, similarly $d(c_1, w_2)=d(v_{1,1}, w_2)$. This allows vertices in $P_1$ and $P_2$ to have distinct vector representations.
     The rest of the proof follows as Theorem \ref{thm:s^m-2,s2,s3}, thus $W$ is a resolving set. 
\end{proof}

\theorem{
    \label{thm:distinct4consecutive}
    Let $G=\Theta(s_1, s_2, s_3, s_4)$, 
    where  $s_1 + 1 = s_2=s_3$ and $s_1+4\leq s_4$. Then $\beta(G)=2$.
}
\begin{proof}
    We claim $W=\{v_{3,1}, v_{4, s_4}\}$ is a resolving set. One can verify that we have that
    \begin{align}
        \vr{c_1} &= (1, s_1 + 2), \\ \vr{c_2} &= (s_1+1, 1),\\
        \vr{v_{1, j}} &= (j + 1, s_1 + 2 - j),\\
        \vr{v_{2, j}} &= (j + 1, s_1 + 3 - j),\\
        \vr{v_{3, j}} &= (j-1, s_1 + 3 - j),
    \end{align}
    from which it is clear that all these vector representations are pairwise distinct. 
    Next, let $u\in V(P_1)\cup V(P_2)$ and $v\in V(P_4)$. Because of our choice of $w_1$, it is always advantageous for $w_1$ to traverse through $c_1$ to reach $u$. It also follows that $w_2$ must not traverse through $c_2$ to reach $v$, otherwise clearly $d(w_2,v)>d(w_2,u)$. Suppose $u$ and $v$ have the same vector representation. Consider \figref{T4}.
    
    \begin{figure}[hbt!]
  \centering
  \begin{tikzpicture}[scale=0.3, transform shape]  
    \input{T4}             
\filldraw[fill=white, draw=black, line width=0.8pt] (v4_2) circle (6pt);
\node[yshift=18pt,  xshift=5pt, scale=2.5] at (v4_2) {$w_1$};

\filldraw[fill=white, draw=black, line width=0.8pt] (v2_2) circle (6pt);
\node[yshift=-18pt,  xshift=10pt, scale=2.5] at (v2_2) {$w_2$};

\filldraw[fill=black, draw=black, line width=0.8pt] (v1_7) circle (6pt);
\node[yshift=18pt,  xshift=10pt, scale=2.5] at (v1_7) {$u$};

\filldraw[fill=black, draw=black, line width=0.8pt] (v2_6) circle (6pt);
\node[yshift=18pt,  xshift=10pt, scale=2.5] at (v2_6) {$v$};
  \end{tikzpicture}
  \caption{\label{fig:T4}}
\end{figure}
Then let $a=d(c_1,u)=d(c_1,v)$ and $b=d(w_2, u)=d(w_2,v)$. Therefore, we get that $s_4=a+b$ and $s_1=a+b-2$, a contradiction as $s_4\geq s_1+4$. Thus $W$ resolves $V(G)\setminus V(P_3)$.

Finally, we have that we may view $C_{3,4}$ independently. By our choice of $w_1$ and $w_2$, both are never MMD to each other and so $W$ resolves $C_{3,4}$. Therefore $W$ is a resolving set.
\end{proof}
\theorem{
    Let $G=\Theta(s_1, s_2, s_3, s_4)$, 
    where  $s_1 + 1 = s_2$ and $s_2<s_3\leq s_4$. Then $\beta(G)=2$.
}
\begin{proof}
    Let $\gamma= \floor{\frac{s_2+s_3+2}{2}}$, then $W=\{w_1,w_2\}$ is a resolving set, where $w_1=v_{3,1}$ and $w_2=v_{3,\gamma}$. We always have that $d(w_2, v_{2,1})= d(w_2,v_{1,1})+1$, which leads to the vertices of $P_1$ having different vector representations to those in $P_2$.

    By our choice of $\gamma$,  $w_1$ and $w_2$ are only MMD when $s_1$ and $s_3$ are of different parities, therefore $C_{1,3}$ is always resolved by $W$.
    A similar argument can be made for $C_{2,3}$.

    Now suppose $u\in V(P_1)\cup V(P_2)$ and $v\in V(P_4)$, such that they have the same vector representation. Both $w_1$ and $w_2$ are MMD to possibly different vertices in $P_4$, say $v_{4,\ell}$ and $v_{4,k}$ respectively. Let $v=v_{4,j}$, then if $j\leq k$ it follows that  
    \begin{align}
        d(w_2, v)>d(w_2, c_1) \geq d(w_2, u),
    \end{align}
    for all $u\in V(P_1)\cup V(P_2)$. A similar argument is made with $w_1$ and $j\geq \ell$. Thus we only need to check when $j\in(k,\ell)$. Consider \figref{T5}.
    \begin{figure}[hbt!]
  \centering
  \begin{tikzpicture}[scale=0.3, transform shape]  
    \input{T4}             
\filldraw[fill=white, draw=black, line width=0.8pt] (v4_2) circle (6pt);
\node[yshift=18pt,  xshift=5pt, scale=2.5] at (v4_2) {$w_1$};

\filldraw[fill=white, draw=black, line width=0.8pt] (v4_7) circle (6pt);
\node[yshift=25pt,  xshift=5pt, scale=2.5] at (v4_7) {$w_2$};

\filldraw[fill=black, draw=black, line width=0.8pt] (v1_7) circle (6pt);
\node[yshift=18pt,  xshift=10pt, scale=2.5] at (v1_7) {$u$};

\filldraw[fill=black, draw=black, line width=0.8pt] (v2_6) circle (6pt);
\node[yshift=18pt,  xshift=10pt, scale=2.5] at (v2_6) {$v$};

\filldraw[fill=black, draw=black, line width=0.8pt] (v2_4) circle (6pt);
\node[yshift=-22pt,  xshift=10pt, scale=2.5] at (v2_4) {$v_{4,\ell}$};

\filldraw[fill=black, draw=black, line width=0.8pt] (v2_8) circle (6pt);
\node[yshift=-24pt,  xshift=-10pt, scale=2.5] at (v2_8) {$v_{4,k}$};
  \end{tikzpicture}
  \caption{\label{fig:T5}}
\end{figure}

    Then $d(w_1,v)=d(w_1,c_1,v)$ and $d(w_2,v)=d(w_2,c_2,v)$.
    Since $w_2$ MMD $c_1$ in $C_{1,3}$, the $d(w_2,u)=d(w_2,c_2,u)$. By our choice of $w_1$, we also have $d(w_1,u)=d(w_1,c_1,u)$. Let $a=d(w_1,u)=d(w_1,v)$ and $b=d(c_2,u)=d(c_2,v)$, then $s_1=a+b-2=s_4$, which is a contradiction. Thus $W$ resolves $V(G)\setminus V(P_3)$.
    
    Lastly we show all vertices in $C_{3,4}$ have unique vector representations. Let $u\in V(P_3)$ and $v\in V(P_4)$, with $u=v_{3,j+1}$ and $v=v_{4,k}$. Suppose $u$ and $v$ have the same vector representation. 
    If $j+1>\gamma$, then since $\gamma> \floor{\frac{s_3}{2}}$, it follows that $d(w_2,v)=d(w_2,c_2,v)$. But $d(w_2,c_2)>d(w_2,u)$,  therefore $j+1<\gamma$. It is therefore apparent that $d(w_2,v)=d(w_2,c_2,v)$. We now use a similar argument as we did for $C_{1,4}$ and $C_{2,4}$. Let $w_1$ MMD $u_1$ and $w_2$ MMD $u_2$ within $P_4$. Then consider \figref{T6}.
        \begin{figure}[hbt!]
  \centering
  \begin{tikzpicture}[scale=0.3, transform shape]  
    \input{T4}             
\filldraw[fill=white, draw=black, line width=0.8pt] (v4_2) circle (6pt);
\node[yshift=18pt,  xshift=5pt, scale=2.5] at (v4_2) {$w_1$};

\filldraw[fill=white, draw=black, line width=0.8pt] (v4_7) circle (6pt);
\node[yshift=25pt,  xshift=5pt, scale=2.5] at (v4_7) {$w_2$};

\filldraw[fill=black, draw=black, line width=0.8pt] (v4_4) circle (6pt);
\node[yshift=18pt,  xshift=10pt, scale=2.5] at (v4_4) {$u$};

\filldraw[fill=black, draw=black, line width=0.8pt] (v2_6) circle (6pt);
\node[yshift=18pt,  xshift=10pt, scale=2.5] at (v2_6) {$v$};

\filldraw[fill=black, draw=black, line width=0.8pt] (v2_4) circle (6pt);
\node[yshift=-22pt,  xshift=10pt, scale=2.5] at (v2_4) {$u_1$};

\filldraw[fill=black, draw=black, line width=0.8pt] (v2_8) circle (6pt);
\node[yshift=-24pt,  xshift=-10pt, scale=2.5] at (v2_8) {$u_2$};
  \end{tikzpicture}
  \caption{\label{fig:T6}}
\end{figure}

    Clearly $d(w_1,v)=d(w_1,c_1,v)$ and as stated $d(w_2,v)=d(w_2,c_2,v)$. Let $a=d(w_1,u)=d(w_1,u)$ and $b=d(w_2,u)=d(w_2,v)$. Then $s_4=a-2 + b - (s_3-\gamma+1)\implies s_4+s_3+3=a+b+\gamma=2\gamma$. But $2\gamma\leq s_2+s_3+2<s_3+s_4+3$, therefore no such $u$ and $v$ exist, and thus $W$ is a resolving set.
\end{proof}
The results in this section are summarized in the following theorem.

\theorem{
    Let $\Theta(s_1, s_2, s_3, s_4)$ be generalized theta graph of multiplicity 4. Then
    \begin{align*}
        \beta(\Theta(s_1, s_2, s_3, s_4)) = \begin{cases}
            4 &\, \text {for $\Theta(1^4)$,$\Theta(1^3, 3)$, $\Theta(2^4)$,or $\Theta(2^3, 4)$}\\
            2 &\, \text{if  $s_2=s_1+1$, $s_2=s_3$ and $s_4\geq s_1 + 4$}
            \\
            2 &\, \text{if  $s_2=s_1+1$ and $s_2<s_3\leq s_4$}
            \\
            3&\,\text{in the remaining cases}
        \end{cases}
    \end{align*}
}

\bibliographystyle{plain} 
\begin{center}
    \bibliography{references} 
\end{center}

\end{document}